\documentclass[11pt]{article}
\usepackage{geometry}                % See geometry.pdf to learn the layout options. There are lots.
\usepackage{graphicx}
\usepackage{amssymb}
\usepackage{amsmath}
\usepackage{epstopdf}
\usepackage[plain]{algorithm2e}
\usepackage{authblk}
\usepackage{tikz}
\newcommand\encircle[1]{%
  \tikz[baseline=(X.base)] 
    \node (X) [draw, shape=circle, inner sep=0] {\strut #1};}
\usepackage[all]{xy}
\title{Two Metrics on Rooted Unordered Trees with Labels}
\author{Yue Wang\thanks{Email address: yuew@g.ucla.edu\\ORCID: 0000-0001-5918-7525}}
\affil{Department of Computational Medicine, University of California, Los Angeles, USA}
\affil{Institut des Hautes \'Etudes Scientifiques, France}

\date{}                                           % Activate to display a given date or no date

\begin{document}
\maketitle

\begin{abstract}
The early development of a zygote can be mathematically described by a developmental tree. To compare developmental trees of different species, we need to define distances on trees. If children cells after a division are not distinguishable, developmental trees are represented by the space $\mathcal{T}$ of rooted trees with possibly repeated labels, where all vertices are unordered. If children cells after a division are partially distinguishable, developmental trees are represented by the space $\mathcal{P}$ of rooted trees with possibly repeated labels, where vertices can be ordered or unordered. On $\mathcal{T}$, the space of rooted unordered trees with possibly repeated labels, we define two metrics: the best-match metric and the left-regular metric, which show some advantages over existing methods. On $\mathcal{P}$, the space of rooted labeled trees with ordered or unordered vertices, there is no metric, and we define a semimetric, which is a variant of the best-match metric. To compute the best-match distance between two trees, the expected time complexity and worst-case time complexity are both $\mathcal{O}(n^2)$, where $n$ is the tree size. To compute the left-regular distance between two trees, the expected time complexity is $\mathcal{O}(n)$, and the worst-case time complexity is $\mathcal{O}(n\log n)$. For rooted labeled trees with (fully/partially) unordered vertices, we define metrics (semimetric) that have fast algorithms to compute and have advantages over existing methods. Such trees also appear outside of developmental biology, and such metrics can be applied to other types of trees which have more extensive applications, especially in molecular biology.
\end{abstract}

\smallskip
\noindent \textbf{Keywords.} 

\noindent Metric; Unordered tree; Label; Semimetric.

\noindent \textbf{MSC code.} 

\noindent 92B05; 05C85.
\section{Background}
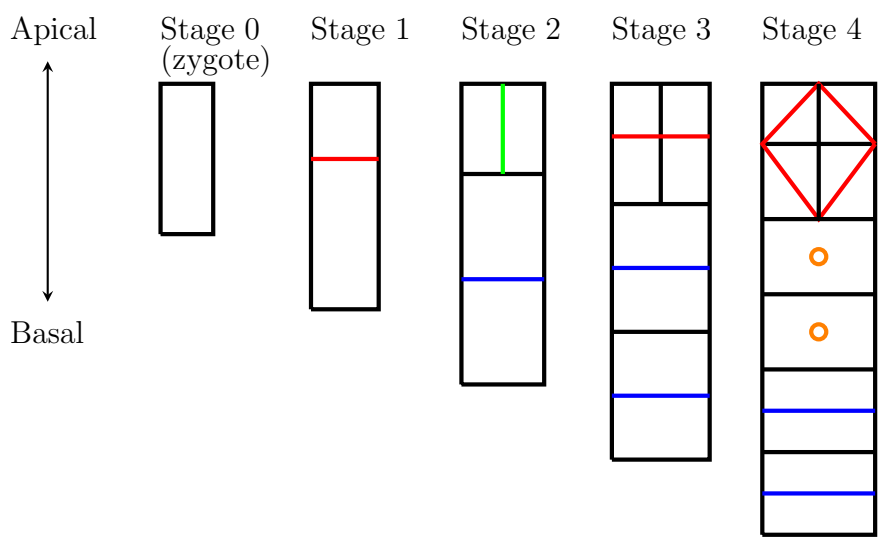
\begin{figure}[h!]
	\begin{tikzpicture}
		\draw [stealth-stealth,thick](-4.5,1.1) -- (-4.5,4.3);
		\node[text width=3cm] at (-3.5,4.7) 
		{\large Apical};
		\node[text width=3cm] at (-3.5,0.7) 
		{\large Basal};
		\node[text width=3cm] at (-1.5,4.7) 
		{\large Stage 0};
		\node[text width=3cm] at (-1.5,4.3) 
		{\large (zygote)};
		\node[text width=3cm] at (0.5,4.7) 
		{\large Stage 1};
		\node[text width=3cm] at (2.5,4.7) 
		{\large Stage 2};
		\node[text width=3cm] at (4.5,4.7) 
		{\large Stage 3};
		\node[text width=3cm] at (6.5,4.7) 
		{\large Stage 4};
		\draw [ultra thick](-3,2) -- (-2.3,2) -- (-2.3,4) -- (-3,4) -- (-3,2);
		\draw [ultra thick](-1,1) -- (-0.1,1) -- (-0.1,4) -- (-1,4) -- (-1,1);
		\draw [ultra thick,red] (-0.1,3) -- (-1,3);
		\draw [ultra thick](1,0) -- (2.1,0) -- (2.1,4) -- (1,4) -- (1,0);
		\draw [ultra thick] (2.1,2.8) -- (1,2.8);
		\draw [ultra thick,blue] (2.1,1.4) -- (1,1.4);
		\draw [ultra thick,green] (1.55,2.8) -- (1.55,4);
		\draw [ultra thick](3,-1) -- (4.3,-1) -- (4.3,4) -- (3,4) -- (3,-1);
		\draw [ultra thick] (4.3,2.4) -- (3,2.4);
		\draw [ultra thick] (4.3,0.7) -- (3,0.7);
		\draw [ultra thick] (3.65,2.4) -- (3.65,4);
		\draw [ultra thick,blue] (4.3,1.55) -- (3,1.55);
		\draw [ultra thick,blue] (4.3,-0.15) -- (3,-0.15);
		\draw [ultra thick,red] (4.3,3.3) -- (3,3.3);
		\draw [ultra thick](5,-2) -- (6.5,-2) -- (6.5,4) -- (5,4) -- (5,-2);
		\draw [ultra thick] (6.5,3.2) -- (5,3.2);
		\draw [ultra thick,red] (6.5,3.2) -- (5.75,4)--(5,3.2)--(5.75,2.2)--(6.5,3.2);
		\draw [ultra thick] (6.5,2.2) -- (5,2.2);
		\draw [ultra thick] (5.75,2.2) -- (5.75,4);
		\draw [ultra thick] (6.5,1.2) -- (5,1.2);
		\draw [ultra thick] (6.5,0.2) -- (5,0.2);
		\draw [ultra thick] (6.5,-0.9) -- (5,-0.9);
		\draw [ultra thick,blue] (6.5,-0.35) -- (5,-0.35);
		\draw [ultra thick,blue] (6.5,-1.45) -- (5,-1.45);
		\draw[color=orange, ultra thick](5.75,0.7) circle (0.1);
		\draw[color=orange, ultra thick](5.75,1.7) circle (0.1);
	\end{tikzpicture}
	\caption{Early development of an \emph{Arabidopsis thaliana} zygote \cite{Wendrich}. Each unit is a cell. A green line between two cells means these two cells were just generated by a symmetric horizontal division. A blue line between two cells means these two cells were just generated by a symmetric vertical division. A red line between two cells means these two cells were just generated by an asymmetric division. An orange circle in a cell means this cell did not divide during the last stage.}
	\label{at}
\end{figure}

In developmental biology, the early development of a zygote is a central topic. For most species, the zygote follows a highly deterministic process. For example, consider a zygote of \emph{Arabidopsis thaliana}. In stage 1, the zygote divides asymmetrically along the apical-basal axis into two cells. In stage 2, the upper (apical) cell undergoes a symmetric horizontal (meridional) division, and the lower (basal) cell undergoes a vertical (equatorial) division. In stage 3, the upper two cells divide asymmetrically, and the lower two cells undergo symmetric vertical divisions. In stage 4, the upper four cells divide asymmetrically, the middle two cells do not divide, and the lower two cells undergo symmetric vertical divisions \cite{Wendrich}. See Fig. \ref{at} for illustrations of this process.

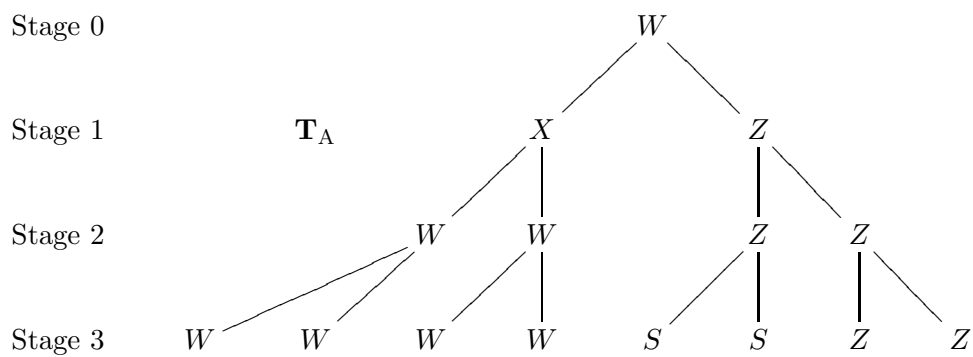
\begin{figure}[t]
	$
	\xymatrix{
		\text{Stage 0}&	&&&&W\ar@{-}[dl]\ar@{-}[dr]&&&\\
		\text{Stage 1}&	&\bf{T_{\rm A}}&&X\ar@{-}[dl]\ar@{-}[d]&&Z\ar@{-}[d]\ar@{-}[dr]&&\\
		\text{Stage 2}&	&&W\ar@{-}[dll]\ar@{-}[dl]&W\ar@{-}[dl]\ar@{-}[d]&&Z\ar@{-}[dl]\ar@{-}[d]&Z\ar@{-}[d]\ar@{-}[dr]&\\
		\text{Stage 3}&	W&W&W&W&S&S&Z&Z
	} 
	$\\
	\caption[]{The developmental tree of \emph{Arabidopsis thaliana}, $T_{\rm A}$, corresponding to Fig. \ref{at}. Each vertex represents a cell, and its label represents the cell event it performs. Label $X$ means symmetric horizontal division; label $Z$ means symmetric vertical division; label $W$ means asymmetric division; label $S$ means the cell stays still and does not divide.}
	\label{atd}
\end{figure}

A mathematical representation of the zygote's early development is a developmental tree \cite{Wang}. In this tree, each vertex represents a cell. Each cell has a label, representing the cell event it will perform, such as division (symmetric or asymmetric, horizontal or vertical), growth, and death. The root vertex is the zygote. Parent vertices (cells) and children vertices (cells) are linked by edges. Each level of this tree corresponds to all the cells at a given stage. See Fig. \ref{atd} for the developmental tree of \emph{Arabidopsis thaliana}.

\begin{figure}[h!]
	\begin{tikzpicture}
		\draw [stealth-stealth,thick](-1.4,-1.6) -- (-1.4,1.6);
		\node[text width=3cm] at (-0.5,2.5) 
		{\large Animal};
		\node[text width=3cm] at (-0.3,2) 
		{\large pole};
		\node[text width=3cm] at (-0.5,-2) 
		{\large Vegetal};
		\node[text width=3cm] at (-0.3,-2.5) 
		{\large pole};
		\node[text width=3cm] at (2,2.4) 
		{\large Stage 0};
		\node[text width=3cm] at (2,-1) 
		{\large Stage 3};
		\node[text width=3cm] at (2,1.9) 
		{\large (zygote)};
		\node[text width=3cm] at (5.6,2.4) 
		{\large Stage 1};
		\node[text width=3cm] at (5.6,-1) 
		{\large Stage 4};
		\node[text width=3cm] at (9.2,2.4) 
		{\large Stage 2};
		\draw [ultra thick](0,0) -- (2.4,0) -- (2.4,1.6) -- (0,1.6) -- (0,0);
		\draw [ultra thick](3.6,0) -- (6,0) -- (6,1.6) -- (3.6,1.6) -- (3.6,0);
		\draw [ultra thick,green] (4.8,0) -- (4.8,1.6);
		\draw [ultra thick] (4.8,-3) -- (4.8,-1.4);
		\draw [ultra thick] (1.2,-3) -- (1.2,-1.4);
		\draw [ultra thick] (4.2,-3) -- (4.2,-1.4);
		\draw [ultra thick] (5.4,-3) -- (5.4,-1.4);
		\draw [ultra thick] (0.6,-3) -- (0.6,-1.4);
		\draw [ultra thick] (1.8,-3) -- (1.8,-1.4);
		\draw [ultra thick] (8.4,0) -- (8.4,1.6);
		\draw [ultra thick,green] (9,0) -- (9,1.6);
		\draw [ultra thick,green] (7.8,0) -- (7.8,1.6);
		\draw [ultra thick](7.2,0) -- (9.6,0) -- (9.6,1.6) -- (7.2,1.6) -- (7.2,0);
		\draw [ultra thick](0,-3) -- (2.4,-3) -- (2.4,-1.4) -- (0,-1.4) -- (0,-3);
		\draw [ultra thick,red](0,-2.1) -- (2.4,-2.1);
		\draw [ultra thick](3.6,-2.1) -- (6,-2.1);
		\draw [ultra thick,red](3.6,-2.7) -- (6,-2.7);
		\draw [ultra thick,green](3.9,-2.1) -- (3.9,-1.4);
		\draw [ultra thick,green](4.5,-2.1) -- (4.5,-1.4);
		\draw [ultra thick,green](5.1,-2.1) -- (5.1,-1.4);
		\draw [ultra thick,green](5.7,-2.1) -- (5.7,-1.4);
		\draw [ultra thick](3.6,-3) -- (6,-3) -- (6,-1.4) -- (3.6,-1.4) -- (3.6,-3);
		
	\end{tikzpicture}
	\caption{Early development of a sea urchin zygote \cite{Summers}. Each unit is a cell. A green line between two cells means these two cells were just generated by a symmetric horizontal division. A red line between two cells means these two cells were just generated by an asymmetric division.}
	\label{su}
\end{figure}
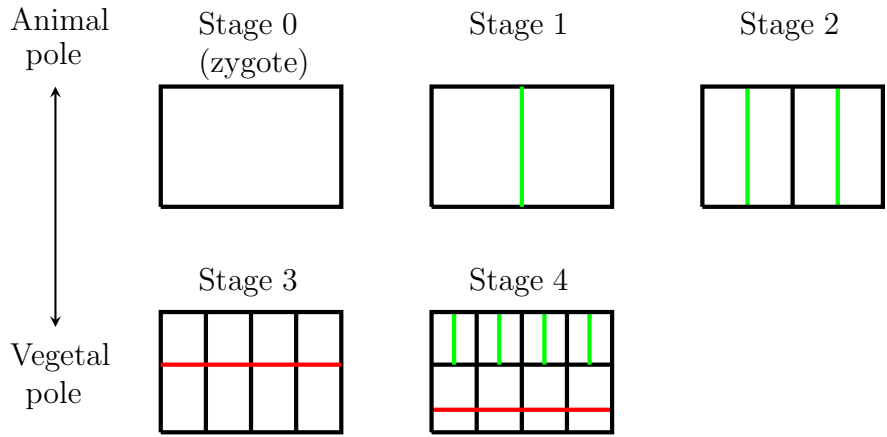

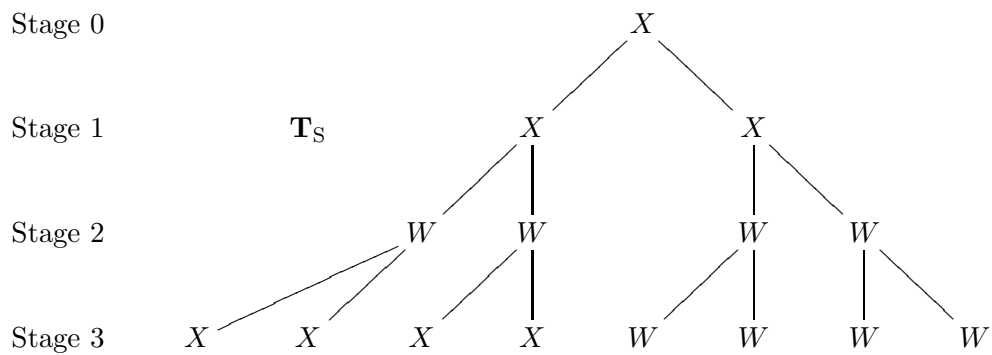
\begin{figure}[t]
	$
	\xymatrix{
		\text{Stage 0}&	&&&&X\ar@{-}[dl]\ar@{-}[dr]&&&\\
		\text{Stage 1}&	&\bf{T_{\rm S}}&&X\ar@{-}[dl]\ar@{-}[d]&&X\ar@{-}[d]\ar@{-}[dr]&&\\
		\text{Stage 2}&	&&W\ar@{-}[dll]\ar@{-}[dl]&W\ar@{-}[dl]\ar@{-}[d]&&W\ar@{-}[dl]\ar@{-}[d]&W\ar@{-}[d]\ar@{-}[dr]&\\
		\text{Stage 3}&	X&X&X&X&W&W&W&W
	} 
	$\\
	\caption[]{The developmental tree of sea urchin, $T_{\rm S}$, corresponding to Fig. \ref{su}. Each vertex represents a cell, and its label represents the cell event it performs. Label $X$ means symmetric horizontal division; label $W$ means asymmetric division.}
	\label{sud}
\end{figure}

Zygotes of different species can have different early developments. See Figs. \ref{su},\ref{sud} for the early development of a sea urchin zygote and the corresponding developmental tree \cite{Summers}. Starting from the zygote, sea urchin and \emph{Arabidopsis thaliana} are different in division plane and division symmetry, and the cell numbers at stage 4 are already different (16 vs. 14). To quantitatively study the development of different organisms, we need a mathematical method to compare different developmental trees. 

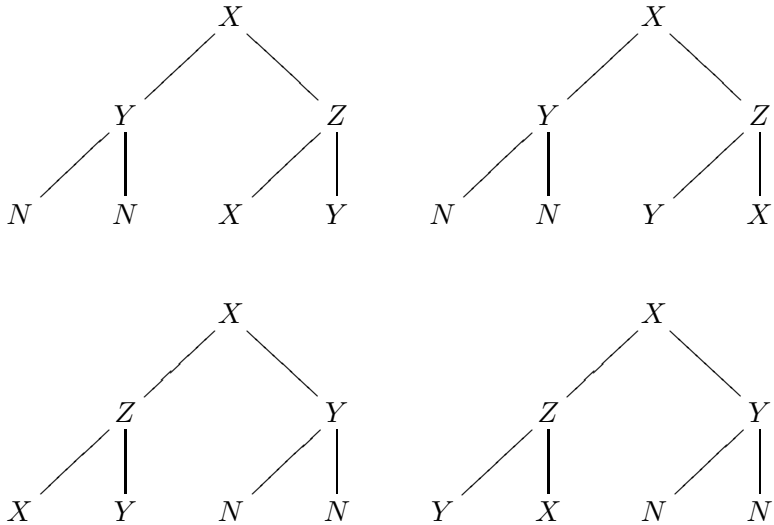
\begin{figure}[t]
	$
	\xymatrix{
		&&X\ar@{-}[dl]\ar@{-}[dr]&&&&X\ar@{-}[dl]\ar@{-}[dr]&\\
		&Y\ar@{-}[dl]\ar@{-}[d]&&Z\ar@{-}[dl]\ar@{-}[d]&&Y\ar@{-}[dl]\ar@{-}[d]&&Z\ar@{-}[dl]\ar@{-}[d]\\
		N&N&X&Y&N&N&Y&X\\
		&&X\ar@{-}[dl]\ar@{-}[dr]&&&&X\ar@{-}[dl]\ar@{-}[dr]&\\
		&Z\ar@{-}[dl]\ar@{-}[d]&&Y\ar@{-}[dl]\ar@{-}[d]&&Z\ar@{-}[dl]\ar@{-}[d]&&Y\ar@{-}[dl]\ar@{-}[d]\\
		X&Y&N&N&Y&X&N&N
	} 
	$\\
	\caption[]{An equivalent class of ordered trees, consisting of four equivalent ordered trees that differ by exchanging the left and right subtrees of some vertices}
	\label{equi}
\end{figure}

When we plot and compare developmental trees, we need to embed them in the plane, namely considering their planar embeddings. We put the zygote to the top, and its two children to the next lower level, and so on. An important question is: after a cell division, which child cell should be put to the left, and which to the right? In some situations, we cannot distinguish two children cells, and we can arbitrarily switch the position of these two children cells in the planar embedding. See Fig. \ref{equi} for equivalent planar embeddings of the same tree. Notice that when we switch cells in the planar embedding, the corresponding cell events are also switched. In some situations, we can distinguish two children cells from an asymmetric division or by which cell inherits the mother centriole \cite{Feldman}. Then we can set a rule to determine which child cell is the left child in the planar embedding, and we cannot switch these two children cells.

We start from the easier situation that we cannot distinguish children cells, so that in the planar embedding of the developmental tree, we can switch two subtrees for each vertex. Notice that a developmental tree has the zygote as its root, and different vertices can have the same label (cell event). The goal is to compare developmental trees.

In the language of graph theory, we need to define a metric on the space of rooted unordered trees with possibly repeated labels. Each tree has a root vertex, and each vertex has a label that is not necessarily unique. All vertices are unordered, meaning that we can switch left and right children in the planar embedding of each tree. Vertices and their labels are always associated, so that we do not distinguish a vertex and the label of a vertex. Therefore, when switching vertices, their labels are also switched. Such trees are not limited to developmental biology, but can be applied in various fields.

There are many metrics defined on trees, which can be roughly classified into three groups by their ideas: (1) Calculate the minimal operations needed to transform one tree into another, such as rearrangement distance \cite{Bernardini}, tree edit distance \cite{Tai}, edge rotation distance \cite{Chartrand}, and geodesic distance \cite{Billera}. (2) Find the largest common structure of two trees, such as bottom-up distance \cite{Valiente} and subtree distance \cite{Zelinka}. (3) Compare structures induced by the trees (e.g., splits or triple-vertices subtrees), such as Robinson-Foulds metric \cite{Robinson}, matching cluster distance \cite{Bogdanowicz}, and triples distance \cite{Critchlow}. 

However, many existing methods have specific requirements on trees, so that they are not applicable in our case (rooted unordered trees with possibly repeated labels). Some methods require that different vertices have different labels, and different trees have the same label set \cite{Bernardini,Chartrand}. Some methods work for phylogenetic trees: only leaves vertices have labels; different vertices have different labels; different trees have the same label set \cite{Billera,Robinson,Bogdanowicz,Critchlow}.
Some methods require that the trees are ordered \cite{Tai}. 

In existing methods, the bottom-up distance \cite{Valiente} and the subtree distance \cite{Zelinka} could work on rooted unordered trees with possibly repeated labels. The bottom-up distance between two trees $T_1,T_2$ is defined as $\mathrm{D_{BU}}(T_1,T_2)=1-f/\max(n_1,n_2)$, where $n_1,n_2$ are the tree sizes, and $f$ is the size of the largest common forest of two trees. The subtree distance $\mathrm{D_{ST}}(T_1,T_2)$ is defined almost the same as the bottom-up distance, except that $f$ is the size of the largest common subtree of two trees. Both distances could be calculated in linear time \cite{flajolet1990analytic,Valiente}. These two methods have some disadvantages. For example, they are not robust under small perturbations on labels, and they do not compare non-common structures. See the next section for detailed discussions. 

We develop two new metrics that apply for rooted unordered trees with possibly repeated labels: the best-match metric $\mathrm{D_{BM}}$ and the left-regular metric $\mathrm{D_{LR}}$. For two unordered trees, the best-match metric searches all their planar embeddings, and compares the most similar pair. To calculate the left-regular metric for two unordered trees, we apply a procedure to fix one planar embedding for each unordered tree (its ``regular form''), and compare the regular forms of these two unordered trees. These two metrics take into account different similarities between labels and different weights concerning their positions. These two metrics, especially the best-match metric, consider any common structures and compare non-common structures. To compute the best-match distance between two trees (binary or general $k$-ary), the expected time complexity and the worst-case time complexity are both $\mathcal{O}(n^2)$, where $n$ is the tree size. To compute the left-regular distance between two trees (binary or not), the expected time complexity is $\mathcal{O}(n)$, and the worst-case time complexity is $\mathcal{O}(n\log n)$. 

The above discussions are for unordered trees, where all vertices are unordered. In some cases, we can distinguish two children cells, so that certain vertices are ordered. Then the space we need to consider consists of rooted trees with possibly repeated labels, where vertices can be ordered or unordered. This larger space has complicated structures that do not allow the existence of a proper metric. Existing methods and the left-regular metric introduced in this paper are not applicable. Nevertheless, the best-match metric can be slightly modified to become a semimetric that works in this scenario.

The main text consists of the following contents: compare existing methods and our new methods; introduce related terminologies in graph theory; define two metrics on the space of rooted unordered trees with possibly repeated labels; define a semimetric on the space of rooted trees with possibly repeated labels, where vertices can be ordered or unordered. 

\section{Comparison of Existing Methods and New Methods}

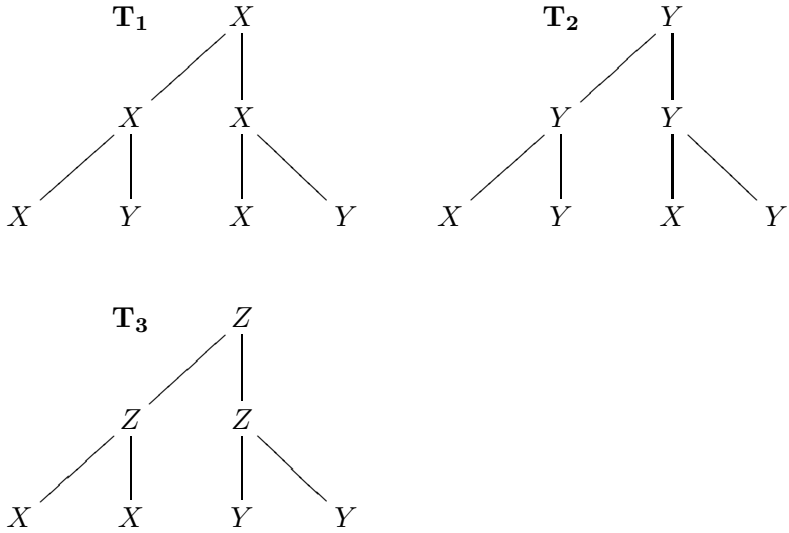
\begin{figure}[t]
	$
	\xymatrix{
		&\bf{T_1}&X\ar@{-}[dl]\ar@{-}[d]&&&\bf{T_2}&Y\ar@{-}[dl]\ar@{-}[d]&\\
		&X\ar@{-}[dl]\ar@{-}[d]&X\ar@{-}[d]\ar@{-}[dr]&&&Y\ar@{-}[dl]\ar@{-}[d]&Y\ar@{-}[d]\ar@{-}[dr]&\\
		X&Y&X&Y&X&Y&X&Y\\
		&\bf{T_3}&Z\ar@{-}[dl]\ar@{-}[d]&&&&&\\
		&Z\ar@{-}[dl]\ar@{-}[d]&Z\ar@{-}[d]\ar@{-}[dr]&&&&&\\
		X&X&Y&Y&&&&
	} 
	$\\
	\caption[]{Three trees $T_1,T_2,T_3$, used to compare the bottom-up distance $\mathrm{D_{BU}}$, the subtree distance $\mathrm{D_{ST}}$, the best-match metric $\mathrm{D_{BM}}$, and the left-regular metric $\mathrm{D_{LR}}$.}
	\label{bu1}
\end{figure}

\begin{figure}[t]
	$
	\xymatrix{
		\bf{T_4}&X\ar@{-}[dl]&&\bf{T_5}&X\ar@{-}[dl]&&\bf{T_6}&X\ar@{-}[dl]\ar@{-}[dr]&\\
		X\ar@{-}[d]&&&Y\ar@{-}[d]&&&Z\ar@{-}[dl]\ar@{-}[d]&&Z\ar@{-}[dl]\ar@{-}[d]\\
		X&&&Y&&Z&Z&Z&Z
	} 
	$\\
	\caption[]{Three trees $T_4,T_5,T_6$, used to compare the bottom-up distance $\mathrm{D_{BU}}$, the subtree distance $\mathrm{D_{ST}}$, the best-match metric $\mathrm{D_{BM}}$, and the left-regular metric $\mathrm{D_{LR}}$.}
	\label{bu2}
\end{figure}
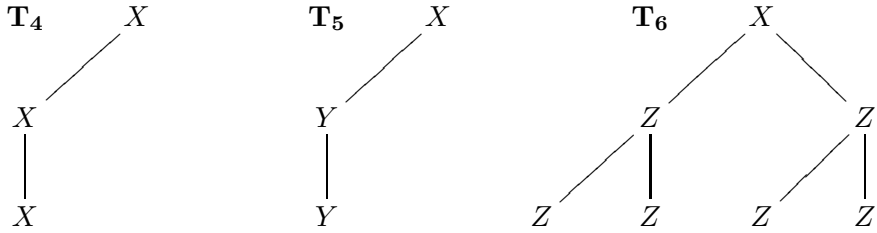

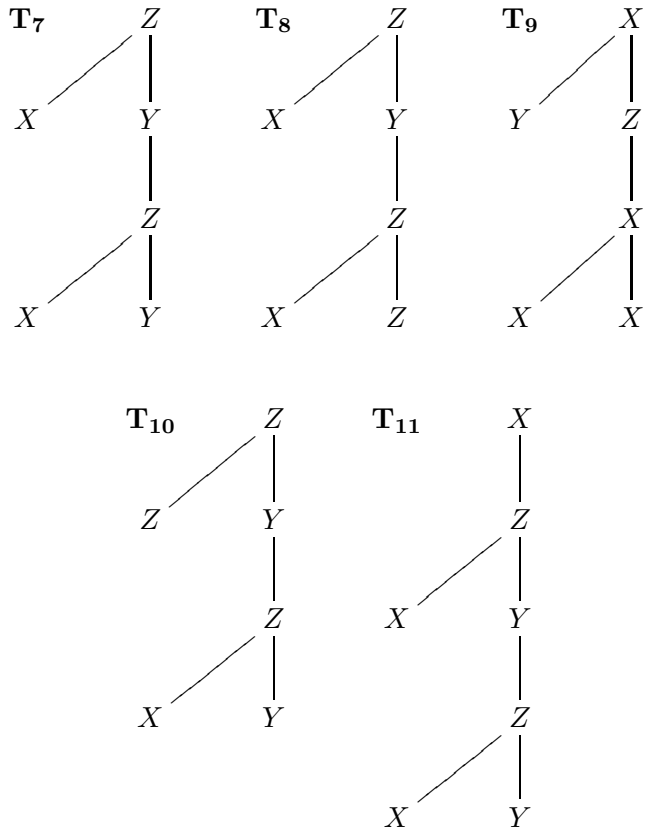
\begin{figure}[t]
	$
	\xymatrix{
		\bf{T_7}&Z\ar@{-}[dl]\ar@{-}[d]&\bf{T_8}&Z\ar@{-}[dl]\ar@{-}[d]&\bf{T_9}&X\ar@{-}[dl]\ar@{-}[d]\\
		X&Y\ar@{-}[d]&X&Y\ar@{-}[d]&Y&Z\ar@{-}[d]\\
		&Z\ar@{-}[dl]\ar@{-}[d]&&Z\ar@{-}[dl]\ar@{-}[d]&&X\ar@{-}[dl]\ar@{-}[d]\\
		X&Y&X&Z&X&X\\
		&\bf{T_{10}}&Z\ar@{-}[dl]\ar@{-}[d]&\bf{T_{11}}&X\ar@{-}[d]&\\
		&Z&Y\ar@{-}[d]&&Z\ar@{-}[dl]\ar@{-}[d]&\\
		&&Z\ar@{-}[dl]\ar@{-}[d]&X&Y\ar@{-}[d]&\\
		&X&Y&&Z\ar@{-}[dl]\ar@{-}[d]&\\
		&&&X&Y&
	} 
	$\\
	\caption[]{Five trees $T_7,T_8,T_9,T_{10},T_{11}$, used to compare the bottom-up distance $\mathrm{D_{BU}}$, the subtree distance $\mathrm{D_{ST}}$, the best-match metric $\mathrm{D_{BM}}$, and the left-regular metric $\mathrm{D_{LR}}$.}
	\label{bu3}
\end{figure}

In this section, we compare the performance of existing methods and new methods on rooted unordered trees with possibly repeated labels, so as to explain the motivation to develop new methods. The examples used are illustrated in Fig. \ref{bu1}, Fig. \ref{bu2}, Fig. \ref{bu3}. See Table \ref{t1} for a summary of these comparisons. 

\begin{table}[]
	\begin{tabular}{lllllll}
		& $\mathrm{D_{BM}}$ & $\mathrm{D_{LR}}$ & \begin{tabular}[c]{@{}l@{}}Size of\\ largest\\ common\\ forest\end{tabular} &$\mathrm{D_{BU}}$  & \begin{tabular}[c]{@{}l@{}}Size of\\ largest\\ common\\ subtree\end{tabular} &$\mathrm{D_{ST}}$  \\
		$T_1$ and $T_2$ & 3 &3  & 4 & 3/7 & 1 &6/7  \\
		$T_1$ and $T_3$ & 5 &  5& 4 & 3/7 &  1& 6/7 \\
		$T_2$ and $T_3$ & 5 &5  & 4 & 3/7 &  1& 6/7  \\
		&&&&&&\\
		$T_4$ and $T_5$ & 2 & 2 &0  & 1 & 0 &1  \\
		$T_4$ and $T_6$ & 6 & 6 & 0 & 1 & 0 & 1 \\
		$T_5$ and $T_6$ & 6 & 6 &0  & 1 & 0 &1  \\
		&&&&&&\\
		$T_7$ and $T_8$ &  1&1  & 2 & 2/3 & 1 &5/6  \\
		$T_7$ and $T_9$ &  5& 5 & 3 & 1/2 &  1&5/6  \\
		$T_7$ and $T_{10}$ & 1 &8  & 4 & 1/3 & 4 & 1/3 \\
		$T_7$ and $T_{11}$ &  9&9  & 6 & 1/7 & 6 & 1/7 \\
		$T_8$ and $T_9$ & 5 & 5 & 2 & 2/3 & 1 &5/6  \\
		$T_8$ and $T_{10}$ & 2 &8  & 2 & 2/3 &1  & 5/6 \\
		$T_8$ and $T_{11}$ & 8 & 8 & 2 & 5/7 &1  &6/7  \\
		$T_9$ and $T_{10}$ &5  &{\color{black}7}  & 2 & 2/3 &1  &5/6  \\
		$T_9$ and $T_{11}$ &7  & 7 &3  &4/7  &1  &6/7  \\
		$T_{10}$ and $T_{11}$ &9  & 10 & 4 &3/7  &4  & 3/7\\
		&&&&&&
	\end{tabular}
	
	\caption[]{Summary of the comparisons in the section ``Comparison of Existing Methods and New Methods''. Performance of $\mathrm{D_{BM}}$, $\mathrm{D_{LR}}$, $\mathrm{D_{BU}}$, and $\mathrm{D_{ST}}$ on trees in Fig. \ref{bu1}, Fig. \ref{bu2}, Fig. \ref{bu3} are illustrated.}
	\label{t1}
\end{table}

Compared to the left-regular metric $\mathrm{D_{LR}}$, especially to the best-match metric $\mathrm{D_{BM}}$ introduced in this paper, the bottom-up distance $\mathrm{D_{BU}}$ \cite{Valiente} and the subtree distance $\mathrm{D_{ST}}$ \cite{Zelinka} have some disadvantages. 

In Fig. \ref{bu1}, $T_1,T_2$ have the same distribution of leaves labels, while $T_1,T_3$ have different distributions of leaves labels. However, $\mathrm{D_{BU}}(T_1,T_2)$$=\mathrm{D_{BU}}(T_1,T_3)=3/7$, $\mathrm{D_{ST}}(T_1,T_2)=\mathrm{D_{ST}}(T_1,T_3)=6/7$. The reason is that $\mathrm{D_{BU}}$ and $\mathrm{D_{ST}}$ only consider \textbf{common structures}, but not their \textbf{detailed patterns}. $\mathrm{D_{BM}}$ and $\mathrm{D_{LR}}$ can recognize the difference: $\mathrm{D_{BM}}(T_1,T_2)=3$, $\mathrm{D_{BM}}(T_1,T_3)=5$; $\mathrm{D_{LR}}(T_1,T_2)=3$, $\mathrm{D_{LR}}(T_1,T_3)=5$. 

In Fig. \ref{bu2}, $T_4,T_5$ have the same {\color{black}tree topology}, while $T_4,T_6$ have different {\color{black}tree topologies}. However, $\mathrm{D_{BU}}(T_4,T_5)=\mathrm{D_{BU}}(T_4,T_6)=1$, $\mathrm{D_{ST}}(T_4,T_6)=\mathrm{D_{ST}}(T_4,T_6)=1$. The reason is that $\mathrm{D_{BU}}$ and $\mathrm{D_{ST}}$ do not compare \textbf{non-common structures}. $\mathrm{D_{BM}}$ and $\mathrm{D_{LR}}$ can recognize the difference: $\mathrm{D_{BM}}(T_4,T_5)=2$, $\mathrm{D_{BM}}(T_4,T_6)=6$; $\mathrm{D_{LR}}(T_4,T_5)=2$, $\mathrm{D_{LR}}(T_4,T_6)=6$. 

In Fig. \ref{bu3}, $T_7,T_8$ only differ by a leaf label, while $T_7,T_9$ are much more different. However, $\mathrm{D_{BU}}(T_7,T_8)=2/3>1/2=\mathrm{D_{BU}}(T_7,T_9)$, $\mathrm{D_{ST}}(T_7,T_8)=\mathrm{D_{ST}}(T_7,T_9)=5/6$. The reason is that $\mathrm{D_{BU}}$ and $\mathrm{D_{ST}}$ only consider \textbf{certain common structures} (sub-forest and subtree). $\mathrm{D_{BM}}$ and $\mathrm{D_{LR}}$ consider any common structures and recognize that $T_7,T_8$ are more similar: $\mathrm{D_{BM}}(T_7,T_8)=1$, $\mathrm{D_{BM}}(T_7,T_9)=5$; $\mathrm{D_{LR}}(T_7,T_8)=1$, $\mathrm{D_{LR}}(T_7,T_9)=5$.

Besides, for two vertices with different labels, $\mathrm{D_{BU}}$ and $\mathrm{D_{ST}}$ only know they are different, but not concerning how different they are. In reality, such as in comparing developmental trees, some labels are very different, while some labels are rather similar. The position of vertices can also be concerned. In general, a label difference closer to the root should be more crucial. In $\mathrm{D_{BM}}$ and $\mathrm{D_{LR}}$, different distances between labels and different weights on vertices can be introduced naturally.

The above discussion explains our motivation to develop the best-match metric $\mathrm{D_{BM}}$ and the left-regular metric $\mathrm{D_{LR}}$. However, $\mathrm{D_{BM}}$ and $\mathrm{D_{LR}}$ also have disadvantages. 

In Fig. \ref{bu3}, $T_7,T_{10}$ only differ by a leaf label. In this case, $\mathrm{D_{BU}}(T_7,T_{10})=1/3$, $\mathrm{D_{ST}}(T_7,T_{10})=1/3$, $\mathrm{D_{BM}}(T_7,T_{10})=1$, but $\mathrm{D_{LR}}(T_7,T_{10})=8$. The reason is that $\mathrm{D_{LR}}$ is not always robust under small perturbations on labels, similar to $\mathrm{D_{BU}}$ and $\mathrm{D_{ST}}$. $\mathrm{D_{BM}}$ is robust under small perturbations on labels. 

In Fig. \ref{bu3}, inserting one vertex to $T_7$ produces $T_{11}$. In this case, $\mathrm{D_{BU}}(T_7,T_{11})=1/7$, $\mathrm{D_{ST}}(T_7,T_{11})=1/7$, but $\mathrm{D_{BM}}(T_7,T_{11})=9$, $\mathrm{D_{LR}}(T_7,T_{11})=9$. The reason is that $\mathrm{D_{BM}}$ and $\mathrm{D_{LR}}$ are not robust under small perturbations on the {\color{black}tree topology}, especially perturbations near the roots. $\mathrm{D_{BU}}$ and $\mathrm{D_{ST}}$ are more robust to the change of {\color{black}tree topology} near the roots.

In summary, our methods outperform the existing methods in most cases. In general, we recommend the best-match metric $\mathrm{D_{BM}}$. If time cost is a major concern, the left-regular metric $\mathrm{D_{LR}}$ can be applied. 

\section{Definitions and Notations}
\subsection{{\color{black}Trees}}
In graph theory, a \textbf{rooted tree} is a connected acyclic undirected graph, where one vertex $v_0$ is designated as the \textbf{root}. Some vertices are linked by edges. For each vertex $v_i$, there is a unique path (edge sequence) that connects $v_i$ and the root $v_0$. The number of edges in this path is called the \textbf{depth} of $v_i$. The depth of the root $v_0$ is stipulated as $0$. The \textbf{depth} of a tree is the largest depth of its vertices. The $k$th \textbf{level} (or level $k$) of a tree consists of all vertices whose depths are $k$. If the depth of a tree is $m$, it is also called an $m$-level tree. If there is an edge between two vertices $v_i,v_j$, and the depth of $v_i$ is smaller than the depth of $v_j$, then $v_i$ is the \textbf{parent} vertex of $v_j$, and $v_j$ is a \textbf{child} vertex of $v_i$. For $v_i$ and its child vertex $v_j$, the tree with root $v_j$ is called a \textbf{subtree} of $v_i$. A vertex without children vertices is called a \textbf{leaf} vertex \cite{Diestel}.

In this paper, each vertex has a \textbf{label}, and different vertices might have the same label. The set of possible labels $\mathcal{L}$ can have infinite elements or even uncountable elements. In the following, we use $\mathcal{L}=\{X,Y,Z\}$ as an example.

For simplicity, we only consider \textbf{binary} trees, meaning that each vertex has at most two children vertices. However, the methods in this paper also work for general $k$-ary trees.

\begin{figure}[t]
	$
	\xymatrix{
		\text{Lv.0}&	\bf{T_{12}}&X\ar@{-}[dl]\ar@{-}[d]&&\bf{\bar{T}_{12}(2)}&&X\ar@{-}[dl]\ar@{-}[dr]&&\\
		\text{Lv.1}&	Y&Z\ar@{-}[d]\ar@{-}[dr]&&&Y\ar@{-}[dl]\ar@{-}[d]&&Z\ar@{-}[d]\ar@{-}[dr]&\\
		\text{Lv.2}&	&Y&Z&N&N&&Y&Z\\
		\text{Lv.0}&	&&&\bf{\bar{T}_{12}(3)}&X\ar@{-}[dl]\ar@{-}[dr]&&&\\
		\text{Lv.1}&	&&&Y\ar@{-}[dl]\ar@{-}[d]&&Z\ar@{-}[d]\ar@{-}[dr]&&\\
		\text{Lv.2}&	&&N\ar@{-}[dll]\ar@{-}[dl]&N\ar@{-}[dl]\ar@{-}[d]&&Y\ar@{-}[dl]\ar@{-}[d]&Z\ar@{-}[d]\ar@{-}[dr]&\\
		\text{Lv.3}&	N&N&N&N&N&N&N&N
	} 
	$\\
	\caption[]{A $2$-level tree $T_{12}$ (upper left), its level-$2$ completion $\bar{T}_{12}(2)$ (upper right) and its level-$3$ completion $\bar{T}_{12}(3)$ (lower). The level (Lv.) of each vertex is marked on the left.}
	\label{com}
\end{figure}
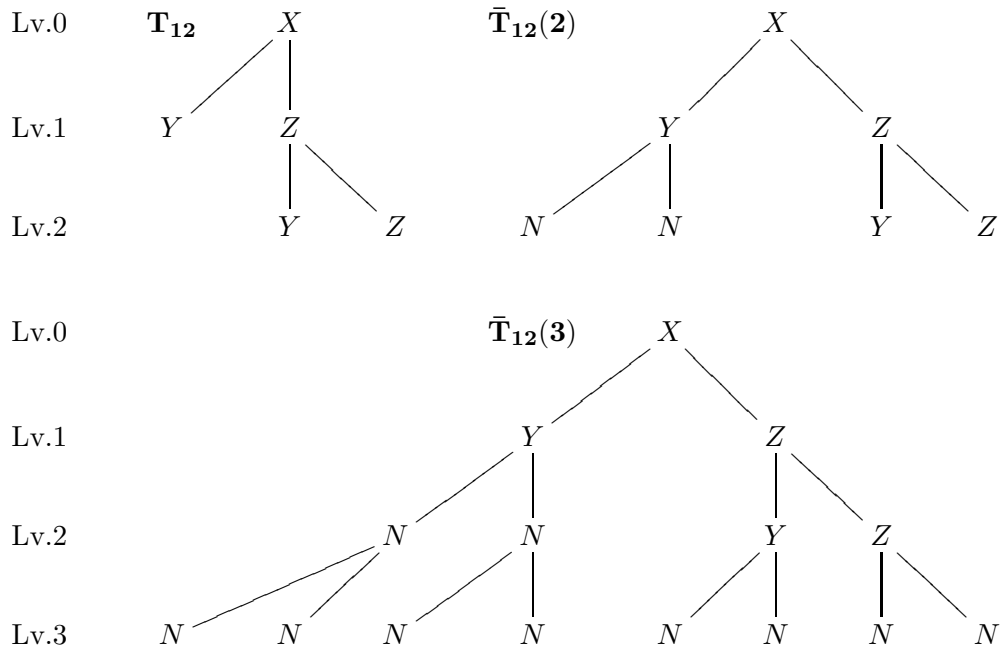

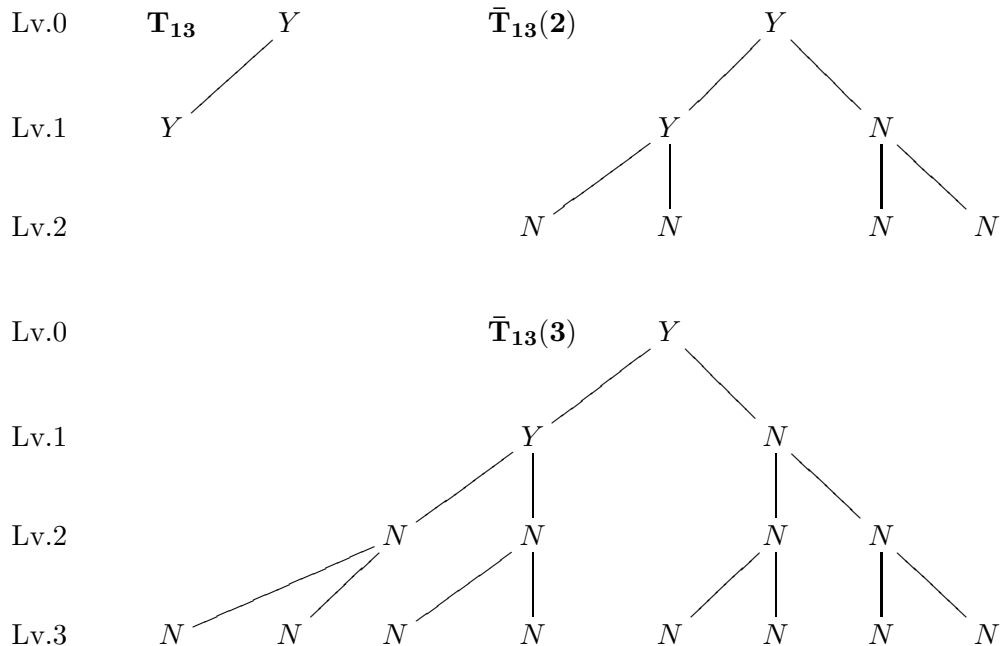
\begin{figure}[t]
	$
	\xymatrix{
		\text{Lv.0}&	\bf{T_{13}}&Y\ar@{-}[dl]&&\bf{\bar{T}_{13}(2)}&&Y\ar@{-}[dl]\ar@{-}[dr]&&\\
		\text{Lv.1}&	Y&&&&Y\ar@{-}[dl]\ar@{-}[d]&&N\ar@{-}[d]\ar@{-}[dr]&\\
		\text{Lv.2}&	&&&N&N&&N&N\\
		\text{Lv.0}&	&&&\bf{\bar{T}_{13}(3)}&Y\ar@{-}[dl]\ar@{-}[dr]&&&\\
		\text{Lv.1}&	&&&Y\ar@{-}[dl]\ar@{-}[d]&&N\ar@{-}[d]\ar@{-}[dr]&&\\
		\text{Lv.2}&	&&N\ar@{-}[dll]\ar@{-}[dl]&N\ar@{-}[dl]\ar@{-}[d]&&N\ar@{-}[dl]\ar@{-}[d]&N\ar@{-}[d]\ar@{-}[dr]&\\
		\text{Lv.3}&	N&N&N&N&N&N&N&N
	} 
	$\\
	\caption[]{A $1$-level tree $T_{13}$ (upper left), its level-$2$ completion $\bar{T}_{13}(2)$ (upper right) and its level-$3$ completion $\bar{T}_{13}(3)$ (lower). The level (Lv.) of each vertex is marked on the left.}
	\label{com2}
\end{figure}

For an $l$-level tree $T$ and any $m\ge l$, we construct its \textbf{level-$\bf \emph m$ completion} $\bar{T}(m)$ as the following: For a vertex not in level $m$, if it has less than two children vertices, add children vertices to it until it has two. Newly added vertices have the label ``$N$'' (means ``null''). Repeat this procedure, until every vertex not in level $m$ has two children vertices, and every vertex in level $m$ has no children vertices. In other words, we construct a perfect binary $m$-level tree. See Fig. \ref{com} and Fig. \ref{com2} for two trees and their completions with different levels. 

For trees after completion, the \textbf{label set} is $\bar{\mathcal{L}}=\mathcal{L}\cup\{N\}$, which is $\{X,Y,Z,N\}$ in our examples. For now, we just require that there is a \textbf{metric} $\mathrm{d}$ on $\bar{\mathcal{L}}$. In this paper, for simplicity, we shall apply the trivial metric that different labels always have distance $1$. Later, we will also need a total order on $\bar{\mathcal{L}}$.

{\color{black}A vertex is called \textbf{ordered} if in the planar embedding of this tree, we know which of its child vertex is the left child, and which is the right child. Otherwise, it is called \textbf{unordered}, and we can switch its two subtrees in the planar embedding. A tree is \textbf{ordered} if all its vertices are ordered. A tree is \textbf{unordered} if all its vertices are unordered.}

Each ordered tree corresponds to a unique planar embedding. In the following, we do not distinguish an ordered tree and its planar embedding. For the space of rooted ordered trees with possibly repeated labels, we define that two trees are \textbf{equivalent} if one tree can transform into the other tree by switching subtrees of some vertices (labels are also switched along with the vertices). Here after transformations, two trees have the same {\color{black}tree topology}, and corresponding vertices have the same label. The notation $T_1\sim T_2$ means $T_1,T_2$ are equivalent, and $T_1\not\sim T_2$ means $T_1,T_2$ are not equivalent. With this equivalence relationship, the space of ordered trees is divided into different equivalent classes. See Fig. \ref{equi} for an equivalent class of ordered trees, where four ordered trees are equivalent. 

An unordered tree corresponds to different planar embeddings (ordered trees). Since we can switch two subtrees of an unordered vertex, equivalent ordered trees represent the same unordered tree. Besides, non-equivalent ordered trees represent different unordered trees. Therefore, the space of unordered trees is isomorphic to the space of equivalent classes of ordered trees. The four ordered trees in Fig. \ref{equi} represent the same unordered tree.

\subsection{{\color{black}Metrics}}
To define a metric on unordered trees, we can switch to equivalent classes of ordered trees. {\color{black}A metric $\mathrm{D}$ on the space of equivalent classes of ordered trees maps a pair of such trees to a non-negative real number, and it satisfies the following \textbf{criteria} for any trees $T_1,T_2,T_3$:

	(A1) $\mathrm{D}(T_1,T_2)=\mathrm{D}(T_2,T_1)$;
	
	(A2) $\mathrm{D}(T_1,T_2)\ge 0$, and $\mathrm{D}(T_1,T_2)=0$ if and only if $T_1\sim T_2$;
	
	(A3) $\mathrm{D}(T_1,T_2)+\mathrm{D}(T_1,T_3)\ge \mathrm{D}(T_2,T_3)$.

	A metric that satisfies (A1)-(A3) also has another property: if $T_1\sim T_2$, then $\mathrm{D}(T_1,T_3)=\mathrm{D}(T_2,T_3)$.}

Before introducing metrics on unordered trees, we first need a metric on the space of ordered trees (not equivalent classes). For two ordered trees $T_1$ and $T_2$, consider their level-$m$ completions, where $m$ is no less than the depths of $T_1$ and $T_2$. For these two completed $m$-level trees $\bar{T}_1(m),\bar{T}_2(m)$ with the same {\color{black}tree topology}, there is a bijection between vertices. We define the \textbf{ordered tree metric} $\mathrm{D_{OT}}(T_1,T_2)$ for such completed ordered trees:
\[\mathrm{D_{OT}}(T_1,T_2)=\mathrm{D_{OT}}(\bar{T}_1(m),\bar{T}_2(m))=\sum_{i\in \bar{T}_1(m)} c(i)\mathrm{d}(i,i'),\]
where $i'\in \bar{T}_2(m)$ is the corresponding vertex of $i$, $\mathrm{d}$ is the metric on the label set $\bar{\mathcal{L}}$, and $c(i)$ is the weight coefficient that depends on the depth of $i$. In some scenarios, we want to emphasize the differences closer to the root (correspond to earlier developmental stages), meaning that we can assign a larger value to $c(i)$ with smaller depth of $i$. For simplicity, we use $c(i)=1$ for all vertices in this paper. We can see that the value of $\mathrm{D_{OT}}$ does not depend on the choice of $m$. For tree $T_{12}$ in Fig. \ref{com} and tree $T_{13}$ in Fig. \ref{com2}, their $\mathrm{D_{OT}}$ distance is 
\[\mathrm{D_{OT}}(T_{12},T_{13})=\mathrm{D_{OT}}(\bar{T}_{12}(2),\bar{T}_{13}(2))=\mathrm{D_{OT}}(\bar{T}_{12}(3),\bar{T}_{13}(3))=4,\]
since they have $4$ pairs of corresponding vertices with different labels. In the rest of this paper, we always consider trees after completion of proper levels. Therefore, the number of vertices (tree size) $n$ and the depth $m$ satisfies $n=2^{m+1}-1$. 

\section{{\color{black}Best-match Metric on Unordered Trees}}
\subsection{{\color{black}Definition}}
We start to define metrics on the space $\mathcal{T}$ of unordered trees, namely the equivalent classes of ordered trees. For two ordered trees $T_1$, $T_2$ (representing their equivalent classes), we can check all pairs of ordered trees that one is equivalent with $T_1$, the other is equivalent with $T_2$, and choose the best-match pair with the minimal $\mathrm{D_{OT}}$ distance. {\color{black}We define $\mathrm{D_{BM}}$ on equivalent classes of ordered trees:
	\[\mathrm{D_{BM}}(T_1,T_2)=\min_{T_1\sim T_1',T_2\sim T_2'}\mathrm{D_{OT}}(T_1',T_2').\]
	This $\mathrm{D_{BM}}(T_1,T_2)$ satisfies the criteria (A1)-(A3) for a metric, defined in the previous section. We name $\mathrm{D_{BM}}$ the \textbf{best-match metric}.} For the tree $T_{12}$ in Fig. \ref{com} and the tree $T_{13}$ in Fig. \ref{com2}, $\mathrm{D_{BM}}(T_{12},T_{13})=4$. 

From the definition of the best-match metric $\mathrm{D_{BM}}(T_1,T_2)$, we can see that changing one label of $T_1$ will make $\mathrm{D_{BM}}(T_1,T_2)$ change by at most $1$. Therefore, the best-match metric is robust under small perturbations on labels. This property does not hold for the left-regular metric, the bottom-up distance, and the subtree distance.

\subsection{{\color{black}A dynamic programming implementation}}
There are exponentially many trees being equivalent to a given tree. Thus brute-force searching is too expensive. Here we introduce a dynamic programming algorithm \cite{ganapathy2006pattern} for calculating the best-match metric $\mathrm{D_{BM}}(T_a,T_b)$.

\begin{algorithm}[!htbp]
	\caption{{\color{black}Detailed workflow of calculating the best-match metric $\mathrm{D_{BM}}$.}}
	\label{bm}
	\vspace{-\bigskipamount}
	%\vspace*{-12pt}
	\ \\
	\begin{enumerate}
		{\color{black}	\item \textbf{Input} 
			
			\quad Two rooted unordered trees $T_a,T_b$ with possibly repeated labels
			
			\quad A metric $\mathrm{d}$ on the label set
			%\quad A function that calculates the ordered tree metric $\mathrm{D_{OT}}$
			
			\item \textbf{Replace} $T_a,T_b$ by the level-$m$ completion of $T_a,T_b$ 
			
			\quad \% Here $m$ is no less than the depths of $T_a$ and $T_b$\
			
			\item \textbf{Define} a function $\mathrm{BM}(l,T_1,T_2)$
			
			\quad \% Here the input is two trees $T_1,T_2$ after completion with depth $l$,
			
			\quad \% and the output is $\mathrm{D_{BM}}(T_1,T_2)$
			
			\quad    \textbf{If} $l=0$
			
			\quad\quad   \textbf{Return} $\mathrm{d}(T_1,T_2)$
			
			\quad   \textbf{Else} 
			
			\quad\quad   \textbf{Define} $T_1^0$ to be $T_1$'s root, and $T_2^0$ to be $T_2$'s root
			
			\quad\quad   \textbf{Define} $T_{1L}$ to be the left subtree of $T_1^0$, 
			
			\quad\quad\quad\quad\quad and $T_{1R}$ to be the right subtree of $T_1^0$
			
			\quad\quad   \textbf{Define} $T_{2L}$ to be the left subtree of $T_2^0$, 
			
			\quad\quad\quad\quad\quad and $T_{2R}$ to be the right subtree of $T_2^0$
			
			\quad\quad   \textbf{Calculate} $C_1=\mathrm{{BM}}(l-1,T_{1L},T_{2L})+\mathrm{{BM}}(l-1,T_{1R},T_{2R})$
			
			\quad\quad   \textbf{Calculate} $C_2=\mathrm{{BM}}(l-1,T_{1L},T_{2R})+\mathrm{{BM}}(l-1,T_{1R},T_{2L})$
			
			\quad\quad   \textbf{Return} $\textrm{d}(T_{1}^0,T_{2}^0)+\min\{C_1,C_2\}$
			
			\quad\quad \% Define $\mathrm{BM}(l,T_1,T_2)$ recursively
			
			\quad   \textbf{End} of if
			
			\item \textbf{Output} $\mathrm{BM}(m,T_a,T_b)$, which is the best-match metric $\mathrm{D_{BM}}(T_a,T_b)$
		}
	\end{enumerate}
	%\vspace{-\bigskipamount}
\end{algorithm}

{\color{black}See Algorithm \ref{bm} for the workflow of calculating the best-match metric $\mathrm{D_{BM}}$. The idea is simple: For the root, we only need to determine whether the left and right subtrees should be switched. In either case, the problem is reduced to minimizing the distance between subtrees. In other words, the vertex correspondence that minimizes the distance between two trees also minimizes the distance between two subtrees. 
	
	In the appendix, we illustrate the detailed procedure of calculating $\mathrm{D_{BM}}$ for the developmental trees of \emph{Arabidopsis thaliana} and sea urchin. $\mathrm{D_{BM}}$ is also applied to other developmental trees with tree size $\sim 100$, and it is discovered that species with similar developmental trees (i.e., smaller $\mathrm{D_{BM}}$) are more likely to have the same anatomical traits \cite{butuzova2022developmental}. For more examples, see Fig. \ref{bu1}, Fig. \ref{bu2}, Fig. \ref{bu3}, and Table \ref{t1}. The Python code for calculating $\mathrm{D_{BM}}$ can be found online (https://github.com/YueWangMathbio/TreeMetric, DOI: 10.5281/zenodo.6400267).} 

\subsection{{\color{black}Computational complexity}}
{\color{black}Assume we need $g(m)$ steps to calculate the best-match distance between two $m$-level trees.} Then we have $g(0)=1$, $g(m+1)=4g(m)+4$. Thus $g(m)=4^{m+1}\times 7/12-4/3$. The number of vertices is $n=2^{m+1}-1$, thus the time complexity of computing the best-match metric is $\mathcal{O}(n^2)$. Here the worst-case time complexity and the expected time complexity are equal. The space complexity of computing the best-match metric is trivially $\mathcal{O}(n)$. When the trees are not binary, but $k$-ary, we have $g(m+1)=k^2g(m)+k\cdot k!\ $. Here the $k^2g(m)$ term means that there are $k^2$ pairs of subtrees to compare. The $k\cdot k!$ term means that for $k!$ possible subtree correspondences, we need $(k-1)\cdot k!$ steps to compute the sum of distances, $k!-1$ steps to compare them, and $1$ step to add $\textrm{d}(T_{1}^0,T_{2}^0)$. Since $g(0)=1$, we have $g(m)=[1/k^2+(k-1)!/(k^2-1)](k^2)^{m+1}-(k\cdot k!)/(k^2-1)$. With $n=(k^{m+1}-1)/(k-1)$, the time complexity is still $\mathcal{O}(n^2)$.

\section{{\color{black}Left-Regular Metric on Unordered Trees}}
\subsection{{\color{black}Preparation}}
Since the metric is defined on the equivalent classes of ordered trees, we need to guarantee that equivalent trees have the same behavior, namely $\mathrm{D}(T_1,T_3)=\mathrm{D}(T_2,T_3)$ for $T_1\sim T_2$. One idea is to transform a given tree into some ``regular form'', which is unique to each equivalent class.

\begin{figure}[t]
	$
	\xymatrix{
		&&&\bf{\bar{T}_{12}(2)}&X\ar@{-}[dl]\ar@{-}[dr]&&&&\\
		&&&Y\ar@{-}[dl]\ar@{-}[d]&&Z\ar@{-}[d]\ar@{-}[dr]&&&\\
		&&N&N&&Y&Z&&\\
		&&X\ar@{-}[dl]\ar@{-}[dr]&&&&&X\ar@{-}[dl]\ar@{-}[dr]&\\
		&Y\ar@{-}[dl]\ar@{-}[d]&&Z\ar@{-}[dr]\ar@{-}[d]&&&Z\ar@{-}[d]\ar@{-}[dl]&&Y\ar@{-}[d]\ar@{-}[dl]\\
		N&N&&Z&Y&Z&Y&N&N
	} 
	$\\
	\caption[]{Procedure of left-regularization. A $2$-level tree $\bar{T}_{12}(2)$ (upper). After left-regularization on level $1$ (lower left). Fully left-regularized (lower right).}
	\label{lreg}
\end{figure}
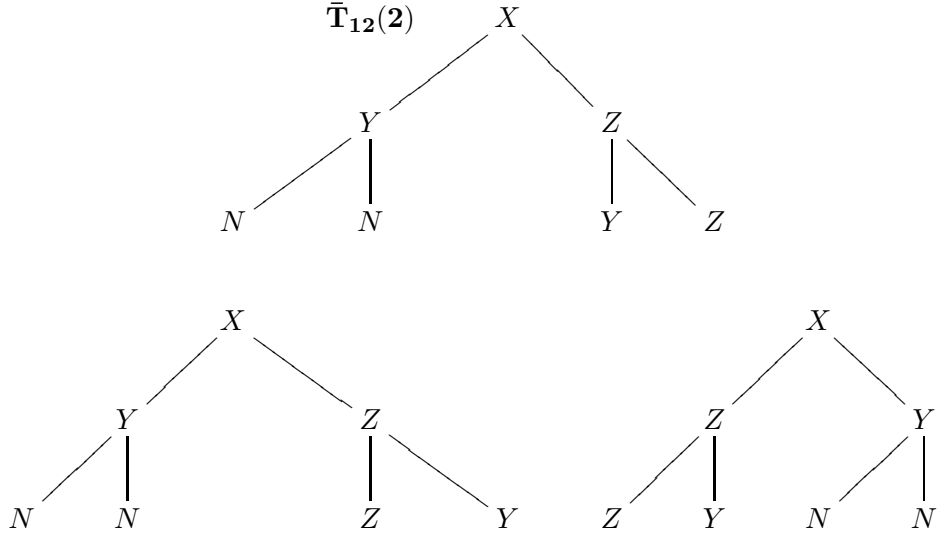

We define a \textbf{total order} on the label set $\bar{\mathcal{L}}$, such as $N>X>Y>Z$. Ideally, similar labels should be closer. With this total order on the label set (alphabet), there is an induced total order, namely the lexicographic order \cite{harzheim2006ordered}, for strings of labels with the same length: for two strings, compare the corresponding labels from the beginning, until there is a difference, and apply the total order for labels. For example, $XZN<XNY$, since $X=X$, and $Z<N$. For a tree after completion, we can write its labels as a string, in the order of up-down (root-leaf), left-right. This is named its \textbf{label string}. For example, the label string of $\bar{T}_{12}(2)$ in Fig. \ref{lreg} is $XYZNNYZ$. We can reconstruct a tree from its label string.

Now we describe the procedure of \textbf{left-regularization}, through which a tree is transformed into its ``regular form''. Consider a tree $T$ after level-$l$ completion. For each vertex in level $l-1$, if the label string of its left subtree is larger than the label string of its right subtree, switch its left and right subtrees. This procedure is called ``left-regularization''. After the left-regularization of level $l-1$, repeat this procedure for level $l-2$,$l-3$,$\ldots$, $1,0$. When the procedure is finished, we obtain the fully ``left-regularized'' form of $T$. The procedure of left-regularization for the tree $\bar{T}_{12}(2)$ is shown in Fig. \ref{lreg}.

In a fully left-regularized tree, for each vertex, the label string of its left subtree is no larger than that of its right subtree. Thus each subtree is also fully left-regularized. By induction with the tree depth, we can see that two equivalent ordered trees have the same left-regularization. Two trees with the same left-regularization are obviously equivalent. Therefore, two ordered trees are equivalent if and only if their left-regularizations are the same. With this procedure, each unordered tree (or its corresponding equivalent class of ordered trees) corresponds to a unique left-regularized ordered tree. 

\subsection{{\color{black}Definition and properties}}
For a tree $T$, denote its fully left-regularized form as $\tilde{T}$. Now we can define $\mathrm{D_{LR}}$ on unordered trees:
\[\mathrm{D_{LR}}(T_1,T_2)=\mathrm{D_{OT}}(\tilde{T}_1,\tilde{T}_2).\]
{\color{black}This $\mathrm{D_{LR}}$ satisfies the criteria (A1)-(A3) for a metric, and we name it the \textbf{left-regular metric}.} Notice that the choice of total order on $\bar{\mathcal{L}}$ might affect the value of $\mathrm{D_{LR}}(T_1,T_2)$. {\color{black}See Algorithm \ref{lr} for the workflow of calculating the left-regular metric $\mathrm{D_{LR}}$. The definition of the left-regular metric already implies how to calculate it. }

\begin{algorithm}[!htbp]
	\caption{{\color{black}Detailed workflow of calculating the left-regular metric $\mathrm{D_{LR}}$.}}
	\label{lr}
	\vspace{-\bigskipamount}
	%\vspace*{-12pt}
	\ \\
	\begin{enumerate}
		{\color{black}	\item \textbf{Input} 
			
			\quad Two rooted unordered trees $T_a,T_b$ with possibly repeated labels
			
			\quad A function that calculates the ordered tree metric $\mathrm{D_{OT}}$
			
			\quad A total order on the label set
			
			\item \textbf{Replace} $T_a,T_b$ by the level-$m$ completion of $T_a,T_b$ 
			
			\quad \% Here $m$ is no less than the depths of $T_a$ and $T_b$\
			
			\item \textbf{For} $l$ from $m-1$ to $0$
			
			\quad \textbf{For} each vertex $v$ in level $l$ of $T_a$, 
			
			\quad \quad \textbf{Set} $S_{\mathrm{L}}(v)$ to be the label string of $v$'s left subtree
			
			\quad \quad \textbf{Set} $S_{\mathrm{R}}(v)$ to be the label string of $v$'s right subtree
			
			\quad \quad   \textbf{If} $S_{\mathrm{L}}(v)>S_{\mathrm{R}}(v)$
			
			\quad\quad\quad   \textbf{Switch} left and right subtrees of $v$
			
			\quad\quad   \textbf{End} of if
			
			\quad   \textbf{End} of for loop
			
			\textbf{End} of for loop
			
			\textbf{Denote} $T_a$ after this left-regularization by $\tilde{T}_a$
			
			\item \textbf{Repeat} Step 3 for $T_b$ to obtain $\tilde{T}_b$
			
			\item \textbf{Output} $\mathrm{D_{OT}}(\tilde{T}_a,\tilde{T}_b)$, which is the left-regular metric $\mathrm{D_{LR}}(T_a,T_b)$
		}
	\end{enumerate}
	%\vspace{-\bigskipamount}
\end{algorithm}

{\color{black}
	The tree $\bar{T}_{13}(2)$ in Fig. \ref{com2} is already left-regularized. Thus we can compare it with $\bar{T}_{12}(2)$ after left-regularization in Fig. \ref{lreg} to find $\mathrm{D_{LR}}(T_{12},T_{13})=5$. In the appendix, we illustrate the detailed procedure of calculating $\mathrm{D_{LR}}$ for the developmental trees of \emph{Arabidopsis thaliana} and sea urchin. For more examples, see Fig. \ref{bu1}, Fig. \ref{bu2}, Fig. \ref{bu3}, and Table \ref{t1}. The Python code for calculating $\mathrm{D_{LR}}$ can be found online (https://github.com /YueWangMathbio/TreeMetric, DOI: 10.5281/zenodo.6400267).} 

Consider an $m$-level tree. For each vertex in level $l$, to compare the label strings of its subtrees, we need at most $2^{m-l}$ steps. Therefore, the left-regularization on each level needs at most $2^m$ steps, and the total number of steps is no more than $(m+1)2^m$. Thus the worst-case time complexity of computing the left-regular metric is $\mathcal{O}(n \log n)$, where $n$ is the vertex number. The space complexity of computing the left-regular metric is trivially $\mathcal{O}(n)$. If the trees are randomly generated, then the expectation of steps needed to compare two label strings is bounded by a constant $C$, regardless of string length. Thus the expected total number of steps is no more than $C2^{m+1}$, and the expected time complexity is $\mathcal{O}(n)$. When the trees are not binary, but $k$-ary, the orders of the worst-case time complexity and the expected time complexity are not changed.

Both the best-match metric and the left-regular metric transform two trees by switching subtrees, and compare the trees after transformation. The best-match metric switches subtrees for two trees cooperatively, so as to find the pair that has the minimal $\mathrm{D}_{\mathrm{OT}}$ distance. The left-regular metric just switches subtrees independently, and the final pair might not be the best match. Thus we can see that {\color{black}for any two unordered trees $T_1,T_2$, $\mathrm{D_{LR}}(T_1,T_2)\ge \mathrm{D_{BM}}(T_1,T_2)$. Thus $\mathrm{D_{LR}}$ is an upper bound of $\mathrm{D_{BM}}$.}

\section{A Semimetric on Trees with Ordered and Unordered Vertices}

\begin{figure}[t]
	$
	\xymatrix{
		\bf{T_{14}^*}&X\ar@{-}[dl]\ar@{-}[dr]&&\bf{T_{15}^*}&\encircle{X}\ar@{-}[dl]\ar@{-}[dr]&&\bf{T_{16}^*}&\encircle{X}\ar@{-}[dl]\ar@{-}[dr]&\\
		Y&&Z\ar@{-}[dl]\ar@{-}[d]&Y&&\encircle{Z}\ar@{-}[dl]\ar@{-}[d]&\encircle{Z}\ar@{-}[dr]\ar@{-}[d]&&Y\\
		&X&Y&&Y&X&X&Y&
	} 
	$\\
	\caption[]{Three ordered trees possibly with lock marks (circles): $T_{14}^*\approx T_{15}^*$, $T_{14}^*\approx T_{16}^*$, $T_{15}^*\not\approx T_{16}^*$.}
	\label{se}
\end{figure}
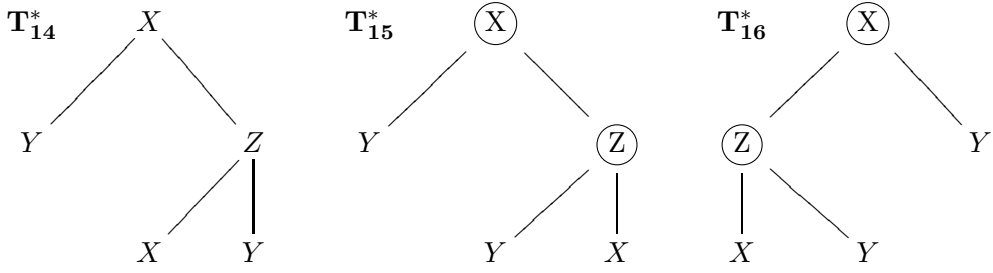

In some situations, for a developmental tree, we know the order of children cells for some cells, but not other cells. Therefore, in this section, we consider the space $\mathcal{P}$ of rooted trees with possibly repeated labels, where vertices can be ordered or unordered. This space contains all ordered trees and unordered trees. To represent this space, we consider ordered trees, where some non-leaves vertices have ``\textbf{lock marks}'', denoted by circles surrounding the labels. See Fig. \ref{se} for some examples. For trees $T_1^*,T_2^*$ possibly with lock marks, if we can switch subtrees of some vertices without lock marks in $T_1^*$, so that $T_1^*$ is transformed into $T_2^*$ (comparing both labels and lock marks), then we define that $T_1^*,T_2^*$ are \textbf{equivalent}, and denote it by $T_1^*\sim T_2^*$. Under this equivalence relation, the space of ordered trees with lock marks is divided into different equivalent classes, and the space of such equivalent classes is isomorphic to the space of trees where vertices can be ordered or unordered. If a tree $T^*$ with lock marks belongs to an equivalent class that represents $T$, then a vertex in $T^*$ has a lock mark if and only if the corresponding vertex in $T$ is ordered. We will define a distance on the equivalent classes of ordered trees with lock marks (not a metric, but a semimetric).

If we try to define a metric on the space of such equivalent classes, we shall meet a problem. Consider two ordered trees $T_1^*,T_2^*$ with the same tree topology and labels, but different lock marks, then $T_1^*\not\sim T_2^*$. If we have a metric $\mathrm{D}$, then $\mathrm{D}(T_1^*,T_2^*)>0$. However, $T_1^*,T_2^*$ have the same labels for corresponding vertices, and we argue that their distance should be $0$. Due to this reason, we need to define another relation between ordered trees with lock marks.

For an ordered tree $T^*$ that could have lock marks, and a normal ordered tree $T$ without lock marks, if we can switch subtrees of some vertices without lock marks in $T^*$, so that $T^*$ is transformed into $T$ (only comparing labels, but not lock marks), then we say $T^*\to T$. For two ordered trees $T_1^*,T_2^*$ that possibly have lock marks, if there exists an ordered tree $T$ without lock marks, so that $T_1^*\to T$ and $T_2^*\to T$, then we define that $T_1^*,T_2^*$ are \textbf{semi-equivalent}, denoted by $T_1^*\approx T_2^*$. $T_1^*\not\approx T_2^*$ means $T_1^*,T_2^*$ are not semi-equivalent. The semi-equivalence relation is reflexive ($T^*\approx T^*$) and symmetric ($T_1^*\approx T_2^*$ means $T_2^*\approx T_1^*$), but not transitive: in Fig. \ref{se}, $T_{14}^*\approx T_{15}^*$, $T_{14}^*\approx T_{16}^*$, but $T_{15}^*\not\approx T_{16}^*$. Thus semi-equivalence is not an equivalence relation, and the space of ordered trees with lock marks cannot be divided into different equivalent classes by this relation. Equivalence is stronger than semi-equivalence: $T_1^*\sim T_2^*$ implies $T_1^*\approx T_2^*$, but not vice versa. Besides, if $T_1^*\sim T_3^*$, $T_2^*\sim T_4^*$, $T_1^*\approx T_2^*$, then $T_3^*\approx T_4^*$.

If $T_1^*\approx T_2^*$, then after certain transformations, they have the same labels, and their distance should be $0$. If $T_1^*\not\approx T_2^*$, then they are essentially different, and their distance should be positive. Therefore, to define a distance $\mathrm{D}$ on the space of ordered trees that could have lock marks, we need to satisfy the following \textbf{criteria} for any trees $T_1^*,T_2^*,T_3^*$: 

(B1) $\mathrm{D}(T_1^*,T_2^*)=\mathrm{D}(T_2^*,T_1^*)$;

{\color{black}(B2) $\mathrm{D}(T_1^*,T_2^*)\ge 0$, and $\mathrm{D}(T_1^*,T_2^*)=0$ if and only if $T_1^*\approx T_2^*$.} 

However, the triangular inequality does not hold. Consider $T_{14}^*,T_{15}^*,T_{16}^*$ in Fig. \ref{se}. Since $T_{14}^*\approx T_{15}^*$, $T_{14}^*\approx T_{16}^*$, $T_{15}^*\not\approx T_{16}^*$, we have $\mathrm{D}(T_{14}^*,T_{15}^*)+\mathrm{D}(T_{14}^*,T_{16}^*)=0<\mathrm{D}(T_{15}^*,T_{16}^*)$. Therefore, we cannot define a metric with respect to the semi-equivalence relation. Besides, $T_1^*\approx T_2^*$ does not imply $\mathrm{D}(T_1^*,T_3^*)=\mathrm{D}(T_2^*,T_3^*)$. Nevertheless, we could require that 

(B3): $T_1^*\sim T_2^*$ implies $\mathrm{D}(T_1^*,T_3^*)=\mathrm{D}(T_2^*,T_3^*)$. 

We seek a distance that satisfies (B1)-(B3), which is a \textbf{semimetric}.

For the space of ordered trees that could have lock marks, existing methods and the left-regular metric are not applicable with respect to the semi-equivalence relation. Nevertheless, the best-match metric can be slightly modified to work in this scenario. {\color{black}On the space of ordered trees that could have lock marks, $\mathrm{D_{BM}^*}$ is defined as
	\[\mathrm{D_{BM}^*}(T_1^*,T_2^*)=\min_{T_1^*\to T_1,T_2^*\to T_2}\mathrm{D_{OT}}(T_1,T_2).\]
	This $\mathrm{D_{BM}^*}(T_1^*,T_2^*)$ satisfies the criteria (B1)-(B3) for a semimetric, and we name it the \textbf{best-match semimetric}.} For the trees in Fig. \ref{se}, $\mathrm{D_{BM}^*}(T_{14}^*,T_{15}^*)=0$, $\mathrm{D_{BM}^*}(T_{14}^*,T_{16}^*)=0$, $\mathrm{D_{BM}^*}(T_{15}^*,T_{16}^*)=6$. {\color{black}See Algorithm \ref{bsm} for the workflow of calculating the best-match semimetric $\mathrm{D_{BM}^*}$. It is only slightly different from Algorithm \ref{bm} that calculates $\mathrm{D_{BM}}$.}

\begin{algorithm}[!htbp]
	\caption{{\color{black}Detailed workflow of calculating the best-match semimetric $\mathrm{D_{BM}^*}$.}}
	\label{bsm}
	\vspace{-\bigskipamount}
	%\vspace*{-12pt}
	\ \\
	\begin{enumerate}
		{\color{black}	\item \textbf{Input} 
			
			\quad Two rooted trees $T_a^*,T_b^*$ with possibly repeated labels,
			
			\quad where some vertices have lock marks
			
			\quad A metric $\mathrm{d}$ on the label set
			%\quad A function that calculates the ordered tree metric $\mathrm{D_{OT}}$
			
			\item \textbf{Replace} $T_a^*,T_b^*$ by the level-$m$ completion of $T_a^*,T_b^*$ 
			
			\quad \% Here $m$ is no less than the depths of $T_a^*$ and $T_b^*$\
			
			\item \textbf{Define} a function $\mathrm{BM}(l,T_1^*,T_2^*)$
			
			\quad \% Here the input is two trees $T_1^*,T_2^*$ after completion with depth $l$,
			
			\quad \% and the output is $\mathrm{D_{BM}^*}(T_1^*,T_2^*)$
			
			\quad    \textbf{If} $l=0$
			
			\quad\quad   \textbf{Return} $\mathrm{d}(T_1^*,T_2^*)$
			
			\quad   \textbf{Else} 
			
			\quad\quad   \textbf{Define} $T_1^{*0}$ to be $T_1^*$'s root, and $T_2^{*0}$ to be $T_2^*$'s root
			
			\quad\quad   \textbf{Define} $T_{1L}^*$ to be the left subtree of $T_1^{*0}$, 
			
			\quad\quad\quad\quad\quad and $T_{1R}^*$ to be the right subtree of $T_1^{*0}$
			
			\quad\quad   \textbf{Define} $T_{2L}^*$ to be the left subtree of $T_2^{*0}$, 
			
			\quad\quad\quad\quad\quad and $T_{2R}^*$ to be the right subtree of $T_2^{*0}$
			
			\quad\quad   \textbf{Calculate} $C_1=\mathrm{{BM}}(l-1,T_{1L}^*,T_{2L}^*)+\mathrm{{BM}}(l-1,T_{1R}^*,T_{2R}^*)$
			
			\quad\quad   \textbf{Calculate} $C_2=\mathrm{{BM}}(l-1,T_{1L}^*,T_{2R}^*)+\mathrm{{BM}}(l-1,T_{1R}^*,T_{2L}^*)$
			
			\quad\quad   \textbf{If} $T_{1}^{*0}$ and $T_{2}^{*0}$ both have lock marks
			
			\quad\quad\quad \textbf{Return} $\textrm{d}(T_{1}^{*0},T_{2}^{*0})+C_1$
			
			\quad\quad   \textbf{Else}
			
			\quad\quad\quad   \textbf{Return} $\textrm{d}(T_{1}^{*0},T_{2}^{*0})+\min\{C_1,C_2\}$
			
			\quad\quad   \textbf{End} of if
			
			\quad\quad \% Define $\mathrm{BM}(l,T_1^*,T_2^*)$ recursively
			
			\quad   \textbf{End} of if
			
			\item \textbf{Output} $\mathrm{BM}(m,T_a^*,T_b^*)$, which is the best-match semimetric $\mathrm{D_{BM}^*}(T_a^*,T_b^*)$
		}
	\end{enumerate}
	%\vspace{-\bigskipamount}
\end{algorithm}

The time complexity of computing the best-match semimetric $\mathrm{D_{BM}^*}$ is also $\mathcal{O}(n^2)$, where $n$ is the number of vertices. The space complexity of computing $\mathrm{D_{BM}^*}$ is $\mathcal{O}(n)$.

\section{Conclusions}

To study the early development of zygotes by comparing developmental trees, we define different distances on trees. On the space of rooted unordered trees with possibly repeated labels, we introduce two metrics: the best-match metric and the left-regular metric. For the same pair of trees, the best-match metric is no larger than the left-regular metric. They consider any common structures and compare non-common structures. Besides, different distances between labels and different weights on vertices can be introduced naturally. The best-match metric has an extra advantage: it is robust under small perturbations on labels. To compute the best-match metric, the time complexity is quadratic, and the left-regular metric has linear expected time complexity. In general, we recommend the best-match metric. If time cost is a major concern, the left-regular metric can be applied. 

On the space of rooted trees with possibly repeated labels, where vertices might be ordered or unordered, most methods are not applicable, and we introduce the best-match semimetric. The properties of the best-match semimetric are almost the same as the best-match metric.

The methods introduced in this paper (possibly with modifications) can be applied to more commonly treated scenarios, e.g., unrooted trees, unlabeled or leaf-labeled trees, or trees with unique labels on vertices. Since our methods are not developed for such scenarios, the performance might not be as satisfactory as existing methods. Nevertheless, the ideas of our methods might inspire new methods in such scenarios.

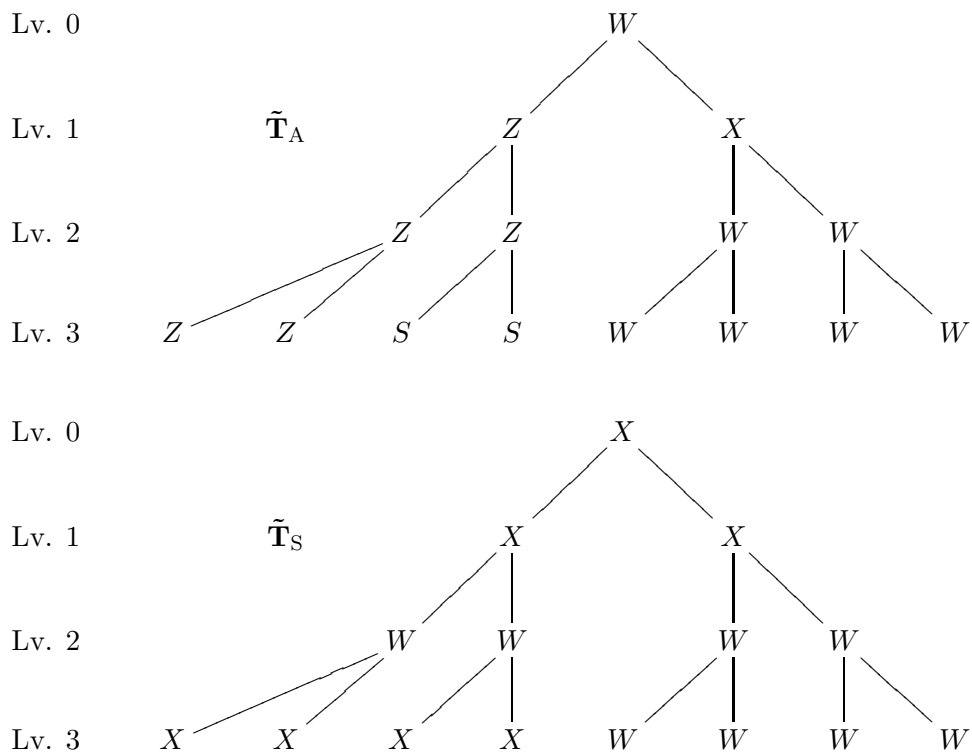
\begin{figure}[t]
	$
	\xymatrix{
		\text{Lv. 0}&	&&&&W\ar@{-}[dl]\ar@{-}[dr]&&&\\
		\text{Lv. 1}&	&\bf{\tilde{T}_{\mathrm{A}}}&&Z\ar@{-}[dl]\ar@{-}[d]&&X\ar@{-}[d]\ar@{-}[dr]&&\\
		\text{Lv. 2}&	&&Z\ar@{-}[dll]\ar@{-}[dl]&Z\ar@{-}[dl]\ar@{-}[d]&&W\ar@{-}[dl]\ar@{-}[d]&W\ar@{-}[d]\ar@{-}[dr]&\\
		\text{Lv. 3}&	Z&Z&S&S&W&W&W&W\\
		\text{Lv. 0}&	&&&&X\ar@{-}[dl]\ar@{-}[dr]&&&\\
		\text{Lv. 1}&	&\bf{\tilde{T}_{\mathrm{S}}}&&X\ar@{-}[dl]\ar@{-}[d]&&X\ar@{-}[d]\ar@{-}[dr]&&\\
		\text{Lv. 2}&	&&W\ar@{-}[dll]\ar@{-}[dl]&W\ar@{-}[dl]\ar@{-}[d]&&W\ar@{-}[dl]\ar@{-}[d]&W\ar@{-}[d]\ar@{-}[dr]&\\
		\text{Lv. 3}&	X&X&X&X&W&W&W&W
	} 
	$\\
	\caption[]{{\color{black}The developmental trees of \emph{Arabidopsis thaliana}, $\tilde{T}_{\mathrm{A}}$ and sea urchin, $\tilde{T}_{\mathrm{S}}$, after left-regularization.}}
	\label{dt}
\end{figure}

{\color{black}\section*{Appendix}
	In this appendix, we present the detailed procedure of calculating the best-match metric $\mathrm{D_{BM}}$ and the left-regular metric $\mathrm{D_{LR}}$ for the developmental trees of \emph{Arabidopsis thaliana} ($T_{\mathrm{A}}$, Fig. \ref{atd}) and sea urchin ($T_{\mathrm{S}}$, Fig. \ref{sud}).
	
	\subsection*{Best-match metric $\mathrm{D_{BM}}$}
	The procedure is recursive. We need to determine the correspondence of subtrees rooted in level 1, which depends on the correspondence of subtrees rooted in level 2, which then depends on the correspondence of subtrees rooted in level 3.
	
	Step 1:
	\begin{equation*}
		\begin{aligned}
			\mathrm{D_{BM}}(T_{\mathrm{A}},T_{\mathrm{S}})=&\mathrm{D_{BM}}(WXZWWZZWWWWSSZZ,\\
			&XXXWWWWXXXXWWWW)\\
			=&\mathrm{d}(W,X)+\min\{\mathrm{D_{BM}}(XWWWWWW,XWWXXXX)\\
			&\qquad\qquad\quad\ \ \,+\mathrm{D_{BM}}(ZZZSSZZ,XWWWWWW),\\
			&\qquad\qquad\qquad\ \ \ \:\! \mathrm{D_{BM}}(XWWWWWW,XWWWWWW)\\
			&\qquad\qquad\quad\ \ \,+\mathrm{D_{BM}}(ZZZSSZZ,XWWXXXX)\}.
		\end{aligned}
	\end{equation*}
	
	Step 2.1:
	\begin{equation*}
		\begin{aligned}
			&\mathrm{D_{BM}}(XWWWWWW,XWWXXXX)\\
			=&\mathrm{d}(X,X)+\min\{\mathrm{D_{BM}}(WWW,WXX)+\mathrm{D_{BM}}(WWW,WXX),\\
			&\qquad\qquad\qquad\quad\!\!\;\mathrm{D_{BM}}(WWW,WXX)+\mathrm{D_{BM}}(WWW,WXX)\}.
		\end{aligned}
	\end{equation*}
	
	Step 3.1.1 to Step 3.1.4 (the same procedure)
	\begin{equation*}
		\begin{aligned}
			&\mathrm{D_{BM}}(WWW,WXX)\\
			=&\mathrm{d}(W,W)+\min\{\mathrm{d}(W,X)+\mathrm{d}(W,X),\mathrm{d}(W,X)+\mathrm{d}(W,X)\}\\
			=&2.
		\end{aligned}
	\end{equation*}
	
	Back to Step 2.1
	\[\mathrm{D_{BM}}(XWWWWWW,XWWXXXX)=4.\]
	
	Step 2.2:
	\begin{equation*}
		\begin{aligned}
			&\mathrm{D_{BM}}(ZZZSSZZ,XWWWWWW)\\
			=&\mathrm{d}(Z,X)+\min\{\mathrm{D_{BM}}(ZSS,WWW)+\mathrm{D_{BM}}(ZZZ,WWW),\\
			&\qquad\qquad\qquad\quad\!\mathrm{D_{BM}}(ZSS,WWW)+\mathrm{D_{BM}}(ZZZ,WWW)\}.
		\end{aligned}
	\end{equation*}
	
	Step 3.2.1 and Step 3.2.3 (the same procedure)
	\begin{equation*}
		\begin{aligned}
			&\mathrm{D_{BM}}(ZSS,WWW)\\
			=&\mathrm{d}(Z,W)+\min\{\mathrm{d}(S,W)+\mathrm{d}(S,W),\mathrm{d}(S,W)+\mathrm{d}(S,W)\}\\
			=&3.
		\end{aligned}
	\end{equation*}
	
	Step 3.2.2 and Step 3.2.4 (the same procedure)
	\begin{equation*}
		\begin{aligned}
			&\mathrm{D_{BM}}(ZZZ,WWW)\\
			=&\mathrm{d}(Z,W)+\min\{\mathrm{d}(Z,W)+\mathrm{d}(Z,W),\mathrm{d}(Z,W)+\mathrm{d}(Z,W)\}\\
			=&3.
		\end{aligned}
	\end{equation*}
	
	Back to Step 2.2
	\[\mathrm{D_{BM}}(ZZZSSZZ,XWWWWWW)=7.\]
	
	Step 2.3:
	\begin{equation*}
		\begin{aligned}
			&\mathrm{D_{BM}}(XWWWWWW,XWWWWWW)\\
			=&\mathrm{d}(X,X)+\min\{\mathrm{D_{BM}}(WWW,WWW)+\mathrm{D_{BM}}(WWW,WWW),\\
			&\qquad\qquad\qquad\quad\!\!\;\mathrm{D_{BM}}(WWW,WWW)+\mathrm{D_{BM}}(WWW,WWW)\}.
		\end{aligned}
	\end{equation*}
	
	Step 3.3.1 to Step 3.3.4 (the same procedure)
	\begin{equation*}
		\begin{aligned}
			&\mathrm{D_{BM}}(WWW,WWW)\\
			=&\mathrm{d}(W,W)+\min\{\mathrm{d}(W,W)+\mathrm{d}(W,W),\mathrm{d}(W,W)+\mathrm{d}(W,W)\}\\
			=&0.
		\end{aligned}
	\end{equation*}
	
	Back to Step 2.3
	\[\mathrm{D_{BM}}(XWWWWWW,XWWWWWW)=0.\]
	
	Step 2.4:
	\begin{equation*}
		\begin{aligned}
			&\mathrm{D_{BM}}(ZZZSSZZ,XWWXXXX)\\
			=&\mathrm{d}(Z,X)+\min\{\mathrm{D_{BM}}(ZSS,WXX)+\mathrm{D_{BM}}(ZZZ,WXX),\\
			&\qquad\qquad\qquad\quad\!\mathrm{D_{BM}}(ZSS,WXX)+\mathrm{D_{BM}}(ZZZ,WXX)\}.
		\end{aligned}
	\end{equation*}
	
	Step 3.4.1 and Step 3.4.3 (the same procedure)
	\begin{equation*}
		\begin{aligned}
			&\mathrm{D_{BM}}(ZSS,WXX)\\
			=&\mathrm{d}(Z,W)+\min\{\mathrm{d}(S,X)+\mathrm{d}(S,X),\mathrm{d}(S,X)+\mathrm{d}(S,X)\}\\
			=&3.
		\end{aligned}
	\end{equation*}
	
	Step 3.4.2 and Step 3.4.4 (the same procedure)
	\begin{equation*}
		\begin{aligned}
			&\mathrm{D_{BM}}(ZZZ,WXX)\\
			=&\mathrm{d}(Z,W)+\min\{\mathrm{d}(Z,X)+\mathrm{d}(Z,X),\mathrm{d}(Z,X)+\mathrm{d}(Z,X)\}\\
			=&3.
		\end{aligned}
	\end{equation*}
	
	Back to Step 2.4
	\[\mathrm{D_{BM}}(ZZZSSZZ,XWWXXXX)=7.\]
	
	Back to Step 1
	\begin{equation*}
		\begin{aligned}
			\mathrm{D_{BM}}(T_{\mathrm{A}},T_{\mathrm{S}})=&\mathrm{D_{BM}}(WXZWWZZWWWWSSZZ,\\
			&XXXWWWWXXXXWWWW)\\
			=&8.
		\end{aligned}
	\end{equation*}

	\subsection*{Left-regularization metric $\mathrm{D_{LR}}$}
	We use the total order $Z<X<W<S$ on the label set. We apply the left-regularization from level $2$ to level $0$.
	
	Left-regularization on level $2$:
	
	For each vertex in level $2$ of $T_{\mathrm{A}}$ and $T_{\mathrm{S}}$, its left subtree and right subtree have the same label string, and we do not need to switch these subtrees. After this step, the label string of $T_{\mathrm{A}}$ is $WXZWWZZWWWWSSZZ$, and the label string of $T_{\mathrm{S}}$ is $XXXWWWWXXXX$ $WWWW$.
	
	Left-regularization on level $1$:
	
	For the vertex with label ``Z'' in level $1$ of $T_{\mathrm{A}}$, its left subtree label string is ``ZSS'', which is larger than that of its right subtree, ``ZZZ''. Thus we switch two subtrees of this vertex. For the other three vertices in level $1$ of $T_{\mathrm{A}}$ and $T_{\mathrm{S}}$, the left subtree and right subtree have the same label string, and we do not need to switch these subtrees. After this step, the label string of $T_{\mathrm{A}}$ is $WXZWWZZWWWWZZSS$, and the label string of $T_{\mathrm{S}}$ is $XXXWWWWXXXXWWWW$.
	
	Left-regularization on level $0$:
	
	For the root vertex (level $0$) of $T_{\mathrm{A}}$, its left subtree label string is ``XWWWWWW'', which is larger than that of its right subtree, ``ZZZZZSS''. Thus we switch two subtrees of this vertex. For the root vertex of $T_{\mathrm{S}}$, its left subtree label string is ``XWWXXXX'', which is smaller than that of its right subtree, ``XWWWWWW''. Thus we do not need to switch these subtrees. After this step, the label string of $T_{\mathrm{A}}$ is $WZXZZWWZZSSWWWW$, and the label string of $T_{\mathrm{S}}$ is $XXXWWWWXXXXWWWW$.
	
	The left-regularization results of $T_{\mathrm{A}}$ and $T_{\mathrm{S}}$ are in Fig. \ref{dt}. We can calculate the $\mathrm{D_{OT}}$ metric for these two trees. Since there are eight pairs of corresponding vertices with different labels, we have $\mathrm{D_{LR}}(T_{\mathrm{A}},T_{\mathrm{S}})=8$.}

\section*{Acknowledgements}
The author would like to thank Mikhail Gromov, Andrey Minarsky, Nadya Morozova, Robert Penner, and Alen Tosenberger for helpful discussions.

\section*{Statements and Declarations}
This research was partially supported by NIH grant R01HL146552.

Declarations of interest: none.

\bibliographystyle{vancouver}
\bibliography{metric}

\end{document}